\documentclass[11pt]{article}
\usepackage[centertags]{amsmath}
\usepackage{amsfonts}
\usepackage{amssymb}
\usepackage{amsthm}
\usepackage{newlfont}
\usepackage{graphicx}

\newfont{\bb}{msbm10 at 12pt}
\def\r{\hbox{\bb R}}

\def\e{\hbox{\bf E}}
\def\t{\hbox{\bf T}}
\def\n{\hbox{\bf N}}
\def\b{\hbox{\bf B}}

\newtheorem{theorem}{Theorem}[section]
\newtheorem{definition}[theorem]{Definition}

\newtheorem{lemma}[theorem]{Lemma}

\begin{document}
\title{\textbf{New special curves and their spherical indicatrices}}

\author{Ahmad T. Ali \\
Mathematics Department,\\
Faculty of Science, Al-Azhar University,\\
Nasr City, 11448, Cairo, Egypt.\\
E-mail: {\tt atali71@yahoo.com}}

\maketitle
\begin{abstract}
In this paper, we define a new special curve in Euclidean 3-space which we call {\it $k-$slant helix} and introduce some characterizations for this curve. This notation is generalization of a general helix and slant helix. Furthermore, we have given some necessary and sufficient conditions for the $k-$slant helix.\\

\textbf{Math. Sub. Class. 2010}: 53A04
\end{abstract}

\section{Introduction}
From the view of differential geometry, a {\it straight line} is a geometric curve with the curvature $\kappa(s)=0$. A {\it plane curve} is a family of geometric curves with torsion $\tau(s)=0$. Helix is a geometric curve with non-vanishing constant curvature $\kappa$ and non-vanishing constant torsion $\tau$ \cite{barros}. The helix may be called a {\it circular helix} or {\it $W$-curve} \cite{ilarslan}. It is known that straight line ($\kappa(s)=0$) and circle ($\kappa(s)=a,\,\tau(s)=0$) are degenerate-helices examples \cite{kuhn}. In fact, circular helix is the simplest three-dimensional spirals \cite{camci}.

A curve of constant slope or {\it general helix} in Euclidean 3-space $\e^3$ is defined by the property that the tangent makes a constant angle with a fixed straight line called the axis of the general helix. A classical result stated by Lancret in 1802 and first proved by de Saint Venant in 1845 (see \cite{struik} for details) says that: {\it A necessary and sufficient condition that a curve be a general helix is that the function
$$f=\dfrac{\tau}{\kappa}$$ is constant along the curve, where $\kappa$ and $\tau$ denote the curvature and the torsion, respectively}. General helices or {\it inclined curves} are well known curves in classical differential geometry of space curves and we refer to the reader for recent works on this type of curves \cite{ali1, ali2, gluck, mont2, turgut}.

In 2004, Izumiya and Takeuchi \cite{izumi} have introduced the concept of {\it slant helix} by saying that the normal lines make a constant angle with a fixed straight line. They characterize a slant helix if and only if the {\it geodesic curvature} of the principal image of the principal normal indicatrix
$$
\sigma=\frac{\kappa^2}{(\kappa^2+\tau^2)^{3/2}}\Big(\frac{\tau}{\kappa}\Big)'
$$
is a constant function. Kula and Yayli \cite{kula1} have studied spherical images of tangent indicatrix and binormal indicatrix of a slant helix and they showed that the spherical images are spherical helices. Recently, Kula et al. \cite{kula2} investigated the relation between a general helix and a slant helix. Moreover, they obtained some differential equations which are characterizations for a space curve to be a slant helix.

A family of curves with constant curvature but non-constant torsion is called Salkowski curves and a family of curves with constant torsion but non-constant curvature is called anti-Salkowski curves $\cite{salkow}$. Monterde \cite{mont1} studied some characterizations of these curves and he proved that the principal normal vector makes a constant angle with fixed straight line. So that: Salkowski and anti-Salkowski curves are the important examples of slant helices.

A unit speed curve of {\it constant precession} in Euclidean 3-space $\e^3$ is defined by the property that its (Frenet) Darboux vector
$$
W=\tau\,\t+\kappa\,\b
$$
revolves about a fixed line in space with constant angle and constant speed. A curve of constant precession is characterized by having
$$
\kappa=\frac{\mu}{m}\sin[\mu\,s],\,\,\,\,\,\,\,\,\,\,\tau=\frac{\mu}{m}\cos[\mu\,s]
$$
or
$$
\kappa=\frac{\mu}{m}\cos[\mu\,s],\,\,\,\,\,\,\,\,\,\,\tau=\frac{\mu}{m}\sin[\mu\,s]
$$
where $\mu$ and $m$ are constants. This curve lies on a circular one-sheeted hyperboloid
$$
x^2+y^2-m^2\,z^2=4m^2
$$
The curve of constant precession is closed if and only if $n=\frac{m}{\sqrt{1+m^2}}$ is rational \cite{scofield}. Kula and Yayli \cite{kula1} proved that the geodesic curvature of the spherical image of the principal normal indicatrix of a curve of constant precession is a constant function equals $-m$. So, one can say that: the curves of constant precessions are the important examples of slant helices.

In this work, we define a new curve and we call it a {\it $k-$slant helix} and we introduce some characterizations of this curve. Furthermore, we have given some necessary and sufficient conditions for the $k-$slant helix. We hope these results will be helpful to mathematicians who are specialized on mathematical modeling as well as other applications of interest.

%%%%%%%%%%%%%%%%%%%%%%%%%%%%%%%%%%%%%%%%%%%%%%%
\section{Preliminaries }
%%%%%%%%%%%%%%%%%%%%%%%%%%%%%%%%%%%%%%%%%%%%%%%
In Euclidean space $\e^3$, it is well known that each unit speed curve with at least four continuous derivatives, one can associate three mutually orthogonal unit vector fields $\t$, $\n$ and $\b$ are respectively, the tangent, the principal normal and the binormal vector fields \cite{hacis}.

We consider the usual metric in Euclidean 3-space $\e^3$, that is,
$$
\langle,\rangle=dx_1^2+dx_2^2+dx_3^2,
$$
where $(x_1,x_2,x_3)$ is a rectangular coordinate system of $\e^3$.  Let $\psi:I\subset\r\rightarrow\e^3$, $\psi=\psi(s)$, be an arbitrary curve in $\e^3$. The curve $\psi$ is said to be of unit speed (or parameterized by the  arc-length) if $\langle\psi'(s),\psi'(s)\rangle=1$ for any $s\in I$. In particular, if $\psi(s)\not=0$ for any $s$, then it is possible to re-parameterize $\psi$, that is, $\alpha=\psi(\phi(s))$ so that $\alpha$ is parameterized by the arc-length. Thus, we will assume throughout this work that $\psi$ is a unit speed curve.

Let $\{\t(s),\n(s),\b(s)\}$ be the moving frame along $\psi$, where the vectors $\t, \n$ and $\b$ are mutually orthogonal vectors satisfying $\langle\t,\t\rangle=\langle\n,\n\rangle=\langle\b,\b\rangle=1$.
The Frenet equations for $\psi$ are given by (\cite{struik,turgut})
\begin{equation}\label{u1}
 \left[
   \begin{array}{c}
     \t'(s) \\
     \n'(s) \\
     \b'(s) \\
   \end{array}
 \right]=\left[
           \begin{array}{ccc}
             0 & \kappa(s) & 0 \\
             -\kappa(s) & 0 & \tau(s) \\
             0 & -\tau(s) & 0 \\
           \end{array}
         \right]\left[
   \begin{array}{c}
     \t(s) \\
     \n(s) \\
     \b(s) \\
   \end{array}
 \right].
 \end{equation}

If $\tau(s)=0$ for all $s\in I$, then $\b(s)$ is a constant vector $V$ and the curve $\psi$ lies in a $2$-dimensional affine subspace orthogonal to $V$, which is isometric to the Euclidean $2$-space $\e^2$.

%%%%%%%%%%%%%%%%%%%%%%%%%%%%%%%%%%%%%%%%%%%%%%%%%%%%%%%%
\section{New representation of spherical indicatrices}
%%%%%%%%%%%%%%%%%%%%%%%%%%%%%%%%%%%%%%%%%%%%%%%%%%%%%%%

In this section we introduce a {\it new representation} of spherical indicatrices of the regular curves in Euclidean 3-space $\e^3$ by the following:

\begin{definition}\label{df-01} Let $\psi$ be a unit speed regular curve in Euclidean 3-space with Frenet vectors $\t$, $\n$ and $\b$. The unit tangent vectors along the curve $\psi(s)$ generate a curve $\psi_{\mathbf{t}}=\t$ on the sphere of radius $1$ about the origin. The curve $\psi_{\mathbf{t}}$ is called the spherical indicatrix of $\t$ or more commonly, $\psi_{\mathbf{t}}$ is called tangent indicatrix of the curve $\psi$. If $\psi=\psi(s)$ is a natural representations of the curve $\psi$, then $\psi_{\mathbf{t}}(s)=\t(s)$ will be a representation of $\psi_{\mathbf{t}}$. Similarly, one can consider the principal normal indicatrix $\psi_{\mathbf{n}}=\n(s)$ and binormal indicatrix $\psi_{\mathbf{b}}=\b(s)$.
\end{definition}

\begin{lemma}\label{lm-01} If the Frenet frame of the tangent indicatrix $\psi_{\mathbf{t}}=\t$ of a space curve $\psi$ is $\{\t_{\mathbf{t}},\n_{\mathbf{t}},\b_{\mathbf{t}}\}$, then we have Frenet formula:
\begin{equation}\label{u2}
 \left[
   \begin{array}{c}
     \t^{\,'}_{\mathbf{t}}(s_{\mathbf{t}})\\
     \n^{\,'}_{\mathbf{t}}(s_{\mathbf{t}})\\
     \b^{\,'}_{\mathbf{t}}(s_{\mathbf{t}})\\
   \end{array}
 \right]=\left[
           \begin{array}{ccc}
             0 & \kappa_{\mathbf{t}} & 0 \\
             -\kappa_{\mathbf{t}} & 0 & \tau_{\mathbf{t}} \\
             0 & -\tau_{\mathbf{t}} & 0 \\
           \end{array}
         \right]\left[
   \begin{array}{c}
     \t_{\mathbf{t}}(s_{\mathbf{t}})\\
     \n_{\mathbf{t}}(s_{\mathbf{t}})\\
     \b_{\mathbf{t}}(s_{\mathbf{t}})\\
   \end{array}
 \right],
 \end{equation}
where
\begin{equation}\label{u3}
\t_{\mathbf{t}}=\n,\,\,\,\,\,\n_{\mathbf{t}}=\frac{-\t+f\,\b}{\sqrt{1+f^2}},\,\,\,\,\,
\b_{\mathbf{t}}=\frac{f\,\t+\b}{\sqrt{1+f^2}},
\end{equation}
and
\begin{equation}\label{u4}
s_{\mathbf{t}}=\int\kappa(s)ds,\,\,\,\,\,\kappa_{\mathbf{t}}=\sqrt{1+f^2},\,\,\,\,\,\tau_{\mathbf{t}}=\sigma\sqrt{1+f^2},
\end{equation}
where
\begin{equation}\label{u41}
f=\frac{\tau(s)}{\kappa(s)}
\end{equation}
and
\begin{equation}\label{u5}
\sigma=\frac{f'(s)}{\kappa(s)\Big(1+f^2(s)\Big)^{3/2}}
\end{equation}
is the geodesic curvature of the principal image of the principal normal indicatrix of the curve $\psi$, $s_{\mathbf{t}}$ is natural representation of the tangent indicatrix of the curve $\psi$ and equal the total curvature of the curve $\psi$ and $\kappa_{\mathbf{t}}$ and $\tau_{\mathbf{t}}$ are the curvature and torsion of $\psi_{\mathbf{t}}$.
\end{lemma}

Therefore we can see that:
\begin{equation}\label{u6}
\frac{\tau_{\mathbf{t}}}{\kappa_{\mathbf{t}}}=\sigma.
\end{equation}

\begin{lemma}\label{lm-02} If the Frenet frame of the principal normal indicatrix $\psi_{\mathbf{n}}=\n$ of a space curve $\psi$ is $\{\t_{\mathbf{n}},\n_{\mathbf{n}},\b_{\mathbf{n}}\}$, then we have Frenet formula:
\begin{equation}\label{u7}
 \left[
   \begin{array}{c}
     \t^{\,'}_{\mathbf{n}}(s_{\mathbf{n}})\\
     \n^{\,'}_{\mathbf{n}}(s_{\mathbf{n}})\\
     \b^{\,'}_{\mathbf{n}}(s_{\mathbf{n}})\\
   \end{array}
 \right]=\left[
           \begin{array}{ccc}
             0 & \kappa_{\mathbf{n}} & 0 \\
             -\kappa_{\mathbf{n}} & 0 & \tau_{\mathbf{n}} \\
             0 & -\tau_{\mathbf{n}} & 0 \\
           \end{array}
         \right]\left[
   \begin{array}{c}
     \t_{\mathbf{n}}(s_{\mathbf{n}})\\
     \n_{\mathbf{n}}(s_{\mathbf{n}})\\
     \b_{\mathbf{n}}(s_{\mathbf{n}})\\
   \end{array}
 \right],
 \end{equation}
 where
\begin{equation}\label{u8}
\left\{
  \begin{array}{ll}
    \t_{\mathbf{n}}=\frac{-\t+f\,\b}{\sqrt{1+f^2}},\\
    \n_{\mathbf{n}}=\frac{\sigma}{\sqrt{1+\sigma^2}}
        \Big[\frac{f\,\t+\b}{\sqrt{1+f^2}}-\frac{\n}{\sigma}\Big],\\
    \b_{\mathbf{n}}=\frac{1}{\sqrt{1+\sigma^2}}\Big[\frac{f\,\t+\b}{\sqrt{1+f^2}}+\sigma\,\n\Big],
  \end{array}
\right.
  \end{equation}
and
\begin{equation}\label{u9}
s_{\mathbf{n}}=\int\kappa(s)\sqrt{1+f^2(s)}\,ds,\,\,\,\,\,\kappa_{\mathbf{n}}=\sqrt{1+\sigma^2},\,\,\,\,\,
\tau_{\mathbf{n}}=\Gamma\sqrt{1+\sigma^2},
\end{equation}
where
\begin{equation}\label{u10}
\Gamma=\frac{\sigma'(s)}{\kappa(s)\sqrt{1+f^2(s)}\Big(1+\sigma^2(s)\Big)^{3/2}},
\end{equation}
$s_{\mathbf{n}}$ is natural representation of the principal normal indicatrix of the curve $\psi$ and $\kappa_{\mathbf{n}}$ and $\tau_{\mathbf{n}}$ are the curvature and torsion of $\psi_{\mathbf{n}}$.
\end{lemma}

Therefore we have:
\begin{equation}\label{u11}
\frac{\tau_{\mathbf{n}}}{\kappa_{\mathbf{n}}}=\Gamma.
\end{equation}

\begin{lemma}\label{lm-03} If the Frenet frame of the binormal indicatrix $\psi_{\mathbf{b}}=\b$ of a space curve $\psi$ is $\{\t_{\mathbf{b}},\n_{\mathbf{b}},\b_{\mathbf{b}}\}$, then we have Frenet formula:
\begin{equation}\label{u12}
 \left[
   \begin{array}{c}
     \t^{\,'}_{\mathbf{b}}(s_{\mathbf{b}})\\
     \n^{\,'}_{\mathbf{b}}(s_{\mathbf{b}})\\
     \b^{\,'}_{\mathbf{b}}(s_{\mathbf{b}})\\
   \end{array}
 \right]=\left[
           \begin{array}{ccc}
             0 & \kappa_{\mathbf{b}} & 0 \\
             -\kappa_{\mathbf{b}} & 0 & \tau_{\mathbf{b}} \\
             0 & -\tau_{\mathbf{b}} & 0 \\
           \end{array}
         \right]\left[
   \begin{array}{c}
     \t_{\mathbf{b}}(s_{\mathbf{b}})\\
     \n_{\mathbf{b}}(s_{\mathbf{b}})\\
     \b_{\mathbf{b}}(s_{\mathbf{b}})\\
   \end{array}
 \right],
 \end{equation}
where
\begin{equation}\label{u13}
\t_{\mathbf{b}}=-\n,\,\,\,\,\,\n_{\mathbf{b}}=\frac{\t-f\,\b}{\sqrt{1+f^2}},\,\,\,\,\,
\b_{\mathbf{b}}=\frac{f\,\t+\b}{\sqrt{1+f^2}},
\end{equation}
and
\begin{equation}\label{u14}
s_{\mathbf{b}}=\int\tau(s)ds,\,\,\,\,\,\kappa_{\mathbf{b}}=\frac{\sqrt{1+f^2}}{f},\,\,\,\,\,
\tau_{\mathbf{b}}=-\frac{\sigma\sqrt{1+f^2}}{f},
\end{equation}
where $s_{\mathbf{b}}$ is natural representation of the binormal indicatrix of the curve $\psi$ and equal the total torsion of the curve $\psi$ and $\kappa_{\mathbf{b}}$ and $\tau_{\mathbf{b}}$ are the curvature and torsion of $\psi_{\mathbf{b}}$.
\end{lemma}

Therefore we obtain:
\begin{equation}\label{u15}
\frac{\tau_{\mathbf{b}}}{\kappa_{\mathbf{b}}}=-\sigma.
\end{equation}

%%%%%%%%%%%%%%%%%%%%%%%%%%%%%%%%%%%%%%%%%%%%%%%%%%%%%%%%
\section{$k$-slant helix and its characterizations}
%%%%%%%%%%%%%%%%%%%%%%%%%%%%%%%%%%%%%%%%%%%%%%%%%%%%%%%
In this section we generalize the concept of the general helix and a slant helix by a new curve which we call it $k$-slant helix.

\begin{definition}\label{df-02} Let $\psi=\psi(s)$ a natural representation of a unit speed regular curve in Euclidean 3-space with Frenet apparatus $\{\kappa,\tau,\t,\n,\b\}$. A curve $\psi$ is called a $k$-slant helix if the unit vector
\begin{equation}\label{u151}
\psi_{\kappa+1}=\frac{\psi'_{k}(s)}{\|\psi'_{k}(s)\|}
\end{equation}
makes a constant angle with a fixed direction, where $\psi_0=\psi(s)$ and $\psi_1=\frac{\psi'_{0}(s)}{\|\psi'_{0}(s)\|}$.
\end{definition}

From the above definition we can see that:

{\bf (1):} The {\it $0$-slant helix} is the curve whose the unit vector
\begin{equation}\label{u16}
\psi_{1}=\frac{\psi'_{0}(s)}{\|\psi'_{0}(s)\|}=\frac{\psi'(s)}{\|\psi'(s)\|}=\t(s),
\end{equation}
(which is the tangent vector of the curve $\psi$) makes a constant angle with a fixed direction. So that the $0$-slant helix is the general helix.

By using the Frenet frame (\ref{u1}), it is easy to prove the following two well-known lemmas:

\begin{lemma}\label{lm-04} Let $\psi:I \rightarrow\e^3$ be a curve that is parameterized by arclength with intrinsic equations $\kappa(s)\neq 0$ and $\tau(s)\neq 0$. The curve $\psi$ is a $0$-slant helix or general helix (the vector $\psi_1$ makes a constant angle, $\phi$, with a fixed straight line in the space) if and only if the function $f(s)=\frac{\tau}{\kappa}=\cot[\phi]$.
\end{lemma}

\begin{lemma}\label{lm-05} Let $\psi:I \rightarrow\e^3$ be a curve that is parameterized by arclength with intrinsic equations $\kappa(s)\neq 0$ and $\tau(s)\neq 0$. The curve $\psi$ is a $0$-slant helix or general helix if and only the binormal vector $\b$ makes a constant angle with fixed direction.
\end{lemma}

{\bf (2):} The {\it $1$-slant helix} is the curve whose the unit vector
\begin{equation}\label{u161}
\psi_{2}=\frac{\psi'_{1}(s)}{\|\psi'_{1}(s)\|}=\frac{\t'(s)}{\|\t'(s)\|}=\n(s),
\end{equation}
(which is the principal normal vector of the curve $\psi$) makes a constant angle with a fixed direction. So that the $1$-slant helix is the slant helix.

If we using the Frenet frame (\ref{u2}) of the tangent indicatrix of the the curve $\psi$, it is easy to prove the following two lemmas. The first lemma is introduced in \cite{ali3, bukcu, izumi, kula1, kula2}. Here, we state this lemma and introduce {\it new representation and its simple proof} using spherical tangent indicatrix of the curve. The second lemma is a new.

\begin{lemma}\label{lm-06} Let $\psi:I\rightarrow\e^3$ be a curve that is parameterized by arclength with intrinsic equations $\kappa(s)\neq0$ and $\tau(s)\neq0$. The curve $\psi$ is a $1$-slant helix or slant helix (the vector $\psi_2$ makes a constant angle, $\phi$, with a fixed straight line in the space) if and only if the function $\sigma(s)=\frac{\tau_{\mathbf{t}}}{\kappa_{\mathbf{t}}}=\cot[\phi]$.
\end{lemma}

{\bf Proof:} $(\Rightarrow)$ Let $\mathbf{d}$ be the unitary fixed vector makes a constant angle, $\phi$, with the vector $\psi_2=\n=\t_{\mathbf{t}}$. Therefore
\begin{equation}\label{u20}
\langle\t_{\mathbf{t}},\mathbf{d}\rangle=\cos[\phi].
\end{equation}
Differentiating the equation (\ref{u20}) with respect to the variable $s_{\mathbf{t}}$ and using Frenet equations (\ref{u2}), we get
\begin{equation}\label{u21}
\kappa_{\mathbf{t}}\langle\n_{\mathbf{t}},\mathbf{d}\rangle=0.
\end{equation}
Because $\kappa_{\mathbf{t}}=\sqrt{1+f^2}\neq0$, then we have
\begin{equation}\label{u22}
\langle\n_{\mathbf{t}},\mathbf{d}\rangle=0.
\end{equation}
From the above equation, the vector $\mathbf{d}$ is perpendicular to the vector $\n_{\mathbf{t}}$ and so that the vector $\mathbf{d}$ lies in the space consists with the vectors $\t_{\mathbf{t}}$ and $\b_{\mathbf{t}}$. Therefore the vector $\mathbf{d}$ makes a constant angles with the two vectors $\t_{\mathbf{t}}$ and $\b_{\mathbf{t}}$. Hence, the vector $\mathbf{d}$ can be written as the following form:
\begin{equation}\label{u23}
\mathbf{d}=\cos[\phi]\t_{\mathbf{t}}+\sin[\phi]\b_{\mathbf{t}}.
\end{equation}
If we differentiate equation (\ref{u23}), we have
\begin{equation}\label{u24}
0=(\cos[\phi]\kappa_{\mathbf{t}}-\sin[\phi]\tau_{\mathbf{t}})\n_{\mathbf{t}},
\end{equation}
which leads to $\sigma(s)=\frac{\tau_{\mathbf{t}}}{\kappa_{\mathbf{t}}}=\cot[\phi]$.

$(\Leftarrow)$ Suppose $\sigma=\cot[\phi]$, i.e., $\tau_{\mathbf{t}}=\cot[\phi]\kappa_{\mathbf{t}}$ and let us consider the vector
\begin{equation}\label{u25}
\mathbf{d}=\cos[\phi]\t_{\mathbf{t}}+\sin[\phi]\b_{\mathbf{t}}.
\end{equation}
From the Frenet formula (\ref{u2}), it is easy to prove the vector $\mathbf{d}$ is constant and $\langle\t_{\mathbf{t}},\mathbf{d}\rangle=\cos[\phi]$. This concludes the proof of lemma (\ref{lm-06}).

\begin{lemma}\label{lm-07} Let $\psi:I \rightarrow\e^3$ be a curve that is parameterized by arclength with intrinsic equations $\kappa(s)\neq 0$ and $\tau(s)\neq 0$. The curve $\psi$ is a $1$-slant helix or slant helix if and only the unit Darboux (modified Darboux \cite{koend}) vector field $\b_{\mathbf{t}}=\frac{f\t+\b}{\sqrt{1+f^2}}$ of $\psi$ makes a constant angle with fixed direction.
\end{lemma}

{\bf Proof:} $(\Rightarrow)$ The proof of the necessary condition is the same as the necessary condition of the above lemma.

$(\Leftarrow)$ Let $\mathbf{d}$ be the unitary fixed vector makes a constant angle, $\frac{\pi}{2}-\phi$, with the vector $\b_{\mathbf{t}}=\frac{f\t+\b}{\sqrt{1+f^2}}$. Therefore
\begin{equation}\label{u26}
\langle\b_{\mathbf{t}},\mathbf{d}\rangle=\sin[\phi].
\end{equation}
Differentiating the equation (\ref{u26}) with respect to the variable $s_{\mathbf{t}}$ and using Frenet equations (\ref{u2}), we get
\begin{equation}\label{u27}
-\tau_{\mathbf{t}}\langle\n_{\mathbf{t}},\mathbf{d}\rangle=0.
\end{equation}
Because $\tau_{\mathbf{t}}=\sigma\sqrt{1+f^2}\neq0$, then we have
\begin{equation}\label{u28}
\langle\n_{\mathbf{t}},\mathbf{d}\rangle=0.
\end{equation}
From the above equation, the vector $\mathbf{d}$ is perpendicular to the vector $\n_{\mathbf{t}}$ and so that the vector $\mathbf{d}$ lies in the space consists with the vectors $\b_{\mathbf{t}}$ and $\t_{\mathbf{t}}$. Therefore the vector $\mathbf{d}$ makes a constant angles with the two vectors $\b_{\mathbf{t}}$ and $\t_{\mathbf{t}}$. This concludes the proof of lemma (\ref{lm-07}).

{\bf (3):} The {\it $2$-slant helix} is the curve whose the unit vector
\begin{equation}\label{u29}
\psi_{3}=\frac{\psi'_{2}(s)}{\|\psi'_{2}(s)\|}=\frac{\n'(s)}{\|\n'(s)\|}=\frac{-\t+f\n}{\sqrt{1+f^2}},
\end{equation}
makes a constant angle with a fixed direction. So that the $2$-slant helix is a new special curves we can call it {\it slant-slant helix}.

If we using the Frenet frame (\ref{u8}) of the principal normal indicatrix of the the curve $\psi$, it is easy to prove the following two new lemmas.

\begin{lemma}\label{lm-08} Let $\psi:I\rightarrow\e^3$ be a curve that is parameterized by arclength with intrinsic equations $\kappa(s)\neq0$ and $\tau(s)\neq0$. The curve $\psi$ is a $2$-slant helix or slant-slant helix (the vector $\psi_3$ makes a constant angle, $\phi$, with a fixed straight line in the space) if and only if the function $\Gamma(s)=\frac{\tau_{\mathbf{n}}}{\kappa_{\mathbf{n}}}=\cot[\phi]$.
\end{lemma}

The proof of the above lemma (using the Frenet frame (\ref{u8})) is similar as the proof of lemma (\ref{lm-06}) (using the Frenet frame (\ref{u2})).

\begin{lemma}\label{lm-09} Let $\psi:I \rightarrow\e^3$ be a curve that is parameterized by arclength with intrinsic equations $\kappa(s)\neq 0$ and $\tau(s)\neq 0$. The curve $\psi$ is a $2$-slant helix or slant-slant helix if and only if the vector $\b_{\mathbf{n}}=\frac{1}{\sqrt{1+\sigma^2}}\Big[\frac{f\,\t+\b}{\sqrt{1+f^2}}+\sigma\,\n\Big]$ makes a constant angle with fixed direction.
\end{lemma}

The proof of the above lemma (using the Frenet frame (\ref{u8})) is similar as the proof of lemma (\ref{lm-07}) (using the Frenet frame (\ref{u2})).

{\bf (4):} The {\it $3$-slant helix} is the curve whose the unit vector
\begin{equation}\label{u30}
\psi_{4}=\frac{\psi'_{3}(s)}{\|\psi'_{3}(s)\|}=\frac{\sigma}{\sqrt{1+\sigma^2}}
        \Big[\frac{f\,\t+\b}{\sqrt{1+f^2}}-\frac{\n}{\sigma}\Big],
\end{equation}
makes a constant angle with a fixed direction. So that the $2$-slant helix is a new special curves we can call it {\it slant-slant-slant helix}.

\begin{lemma}\label{lm-10} Let $\psi:I\rightarrow\e^3$ be a curve that is parameterized by arclength with intrinsic equations $\kappa(s)\neq0$ and $\tau(s)\neq0$. The curve $\psi$ is a $3$-slant helix or slant-slant-slant helix (the vector $\psi_4$ makes a constant angle, $\phi$, with a fixed straight line in the space) if and only if the function \begin{equation}\label{u301}
\Lambda=\frac{\Gamma'(s)}{\kappa(s)\sqrt{1+f^2(s)}\sqrt{1+\sigma^2(s)}\Big(1+\Gamma^2(s)\Big)^{3/2}}=\cot[\phi].
\end{equation}
\end{lemma}

{\bf proof:} $(\Rightarrow)$ Let $\textbf{d}$ be the unitary fixed vector makes a constant angle, $\phi$, with the vector $\psi_{4}=\n_{\mathbf{n}}$. Therefore
\begin{equation}\label{u31}
\langle\n_{\mathbf{n}},\textbf{d}\rangle=\cos[\phi].
\end{equation}
Differentiating the equation (\ref{u31}) with respect to the variable $s_{\mathbf{n}}=\int\kappa(s)\sqrt{1+f^2(s)}ds$ and using the Frenet equations (\ref{u8}), we get
\begin{equation}\label{u32}
\langle-\kappa_{\mathbf{n}}\t_{\mathbf{n}}+\tau_{\mathbf{n}}\b_{\mathbf{n}},\textbf{d}\rangle=0.
\end{equation}
Therefore,
$$
\langle\t_{\mathbf{n}},\textbf{d}\rangle=\frac{\tau_{\mathbf{n}}}{\kappa_{n}}\langle\b_{\mathbf{n}},\textbf{d}\rangle=
\Gamma\langle\b_{\mathbf{n}},\textbf{d}\rangle.
$$
If we put $\langle\b_{\mathbf{n}},\textbf{d}\rangle=g(s)$, we can write
$$
\textbf{d}=\Gamma\,g\,\t_{\mathbf{n}}+\cos[\phi]\n_{\mathbf{n}}+g\,\b_{n}.
$$
From the unitary of the vector $\textbf{d}$ we get $g=\pm
\frac{\sin[\phi]}{\sqrt{1+\Gamma^2}}$. Therefore, the vector $\textbf{d}$ can be written as
\begin{equation}\label{u33}
\textbf{d}=\pm\,\frac{\Gamma\,\sin[\phi]}{\sqrt{1+\Gamma^2}}\,\t_{\mathbf{n}}+\cos[\phi]\,\n_{\mathbf{n}}
\pm\frac{\sin[\phi]}{\sqrt{1+\Gamma^2}}\,\b_{\mathbf{n}}.
\end{equation}
The equation (\ref{u32}) can be written in the form:
\begin{equation}\label{u34}
\langle-\t_{\mathbf{n}}+\Gamma\,\b_{\mathbf{n}},\textbf{d}\rangle=0.
\end{equation}
If we differentiate the equation (\ref{u32}) with respect to $s_{\mathbf{n}}$, again, we obtain
\begin{equation}\label{u35}
\langle \dot{\Gamma}\,\b_{\mathbf{n}}+(1+\Gamma^2)\sqrt{1+\sigma^2}\n_{\mathbf{n}},\textbf{d}\rangle=0,
\end{equation}
where dot is the differentiation with respect to $s_{\mathbf{n}}$. If we put the vector $\mathbf{d}$ from equation (\ref{u33}) in the equation (\ref{u35}), we obtain the following condition
$$
\frac{\dot{\Gamma}}{\sqrt{1+\sigma^2}(1+\Gamma^2)^{3/2}}=\pm\,\cot[\phi].
$$
Finally, $s_{\mathbf{n}}=\int\kappa(s)\sqrt{1+f^2(s)}ds$ and $\dot{\Gamma}=\frac{\Gamma'(s)}{\kappa(s)\sqrt{1+f^2(s)}}$, we express the desired result.

$(\Leftarrow)$ Suppose that $\frac{\dot{\Gamma}}{\sqrt{1+\sigma^2}(1+\Gamma^2)^{3/2}}=\pm\,\cot[\phi]$ where $.$ is the differentiation with respect to $s_{\mathbf{n}}$. Let us consider the vector
$$
\textbf{d}=\pm\,\cos[\phi]\Big(\frac{\Gamma\,\tan[\phi]}{\sqrt{1+\Gamma^2}}\,\t_{\mathbf{n}}\pm\n_{\mathbf{n}}
+\frac{\tan[\phi]}{\sqrt{1+\Gamma^2}}\,\b_{\mathbf{n}}\Big).
$$
We will prove that the vector $\textbf{d}$ is a constant vector. Indeed, applying Frenet formula (\ref{u8})
$$
\dot{\textbf{d}}=\pm\sqrt{1+\sigma^2}\cos[\phi]\Big(\pm\t_{\mathbf{n}}+\frac{\Gamma\tan[\phi]}{\sqrt{1+\Gamma^2}}\,\n_{n}
\mp\t_{\mathbf{n}}\pm\Gamma\b\mp\Gamma\b_{\mathbf{n}}
-\frac{\Gamma\tan[\phi]}{\sqrt{1+\Gamma^2}}\,\n_{n}\Big)=0
$$
Therefore, the vector $\textbf{d}$ is constant and $\langle\n_{\mathbf{n}},\textbf{d}\rangle=\cos[\phi]$. This concludes the proof of lemma (\ref{lm-10}).

From the section (3), we can see that:

{\bf (i):} The function $f(s)$ is equal the ratio of the torsion $(\tau=\tau_0)$ and curvature $(\kappa=\kappa_0)$ of the curve $\psi=\psi_0$ and may be named it $\sigma_0(s)=f(s)=\frac{\tau_0(s)}{\kappa_0(s)}$.

{\bf (ii):} The function $\sigma(s)$ is equal the ratio of the torsion $(\tau_{\mathbf{t}}=\tau_1)$ and curvature $(\kappa_{\mathbf{t}}=\kappa_1)$ of the tangent indicatrix $\t=\psi_1$ of the curve $\psi$ and may be named it $\sigma_1(s)=\sigma(s)=\frac{\tau_1(s)}{\kappa_1(s)}$.

{\bf (iii):} The function $\Gamma(s)$ is equal the ratio of the torsion $(\tau_{\mathbf{n}}=\tau_2)$ and curvature $(\kappa_{\mathbf{n}}=\kappa_2)$ of the principal normal indicatrix $\n=\psi_2$ of the curve $\psi$ and may be named it $\sigma_2(s)=\Gamma(s)=\frac{\tau_2(s)}{\kappa_2(s)}$.

We expect that: the function $\Lambda(s)$ is equal the ratio of the torsion $\tau_3$ and curvature $\kappa_3$ of the spherical image of $\psi_3$ indicatrix and may be named it $\sigma_3(s)=\Lambda(s)=\frac{\tau_3(s)}{\kappa_3(s)}$. So that, we can write (the proof is classical) the following lemma:

\begin{lemma}\label{lm-11} If the Frenet frame of the spherical image of $\psi_3=\frac{-\t+f\b}{\sqrt{1+f^2}}$ indicatrix of the curve $\psi$ is $\{\t_3,\n_3,\b_3\}$, then we have Frenet formula:
\begin{equation}\label{u36}
 \left[
   \begin{array}{c}
     \t^{\,'}_3(s_3)\\
     \n^{\,'}_3(s_3)\\
     \b^{\,'}_3(s_3)\\
   \end{array}
 \right]=\left[
           \begin{array}{ccc}
             0 & \kappa_3 & 0 \\
             -\kappa_3 & 0 & \tau_3 \\
             0 & -\tau_3 & 0 \\
           \end{array}
         \right]\left[
   \begin{array}{c}
     \t_3(s_3)\\
     \n_3(s_3)\\
     \b_3(s_3)\\
   \end{array}
 \right],
 \end{equation}
 where
\begin{equation}\label{u37}
\left\{
  \begin{array}{ll}
        \t_3=\frac{\sigma}{\sqrt{1+\sigma^2}}
        \Big[\frac{f\,\t+\b}{\sqrt{1+f^2}}-\frac{\n}{\sigma}\Big],\\
    \n_3=\frac{1}{\sqrt{1+\sigma^2}\sqrt{1+\Gamma^2}}\Big[
    \frac{\Gamma\big(f\t+\b\big)+\sqrt{1+\sigma^2}\big(\t-f\b\big)}{\sqrt{1+f^2}}+\sigma\Gamma\n\Big],\\
    \b_3=\frac{1}{\sqrt{1+\sigma^2}\sqrt{1+\Gamma^2}}\Big[
    \frac{f\t+\b-\Gamma\sqrt{1+\sigma^2}\big(\t-f\b\big)}{\sqrt{1+f^2}}+\sigma\n\Big],
  \end{array}
\right.
  \end{equation}
and
\begin{equation}\label{u38}
s_3=\int\kappa(s)\sqrt{1+f^2(s)}\sqrt{1+\sigma^2(s)}\,ds,\,\,\,\,\,\kappa_3=\sqrt{1+\Gamma^2},\,\,\,\,\,
\tau_3=\Lambda\sqrt{1+\Gamma^2},
\end{equation}
where
$s_3$ is the natural representation of the spherical image of $\psi_3$ indicatrix of the curve $\psi$ and $\kappa_3$ and $\tau_3$ are the curvature and torsion of this curve.
\end{lemma}

Therefore it is easy to see that:
\begin{equation}\label{u39}
\frac{\tau_3}{\kappa_3}=\Lambda=\sigma_3.
\end{equation}

If we using the Frenet frame (\ref{u37}) of the spherical image of $\psi_3$ indicatrix of the curve $\psi$, it is easy to prove the following new lemma.

\begin{lemma}\label{lm-12} Let $\psi:I \rightarrow\e^3$ be a curve that is parameterized by arclength with intrinsic equations $\kappa(s)\neq 0$ and $\tau(s)\neq 0$. The curve $\psi$ is a $3$-slant helix or slant-slant-slant helix if and only if the vector $\b_3=\frac{1}{\sqrt{1+\sigma^2}\sqrt{1+\Gamma^2}}\Big[
    \frac{f\t+\b-\Gamma\sqrt{1+\sigma^2}\big(\t-f\b\big)}{\sqrt{1+f^2}}+\sigma\n\Big]$ makes a constant angle with fixed direction.
\end{lemma}

The proof of the above lemma (using the Frenet frame (\ref{u37})) is similar as the proof of lemma (\ref{lm-07}) (using the Frenet frame (\ref{u2})).

%%%%%%%%%%%%%%%%%%%%%%%%%%%%%%%%%%%%%%%%%%%%%%%%%%%%%%%%
\section{General results}
%%%%%%%%%%%%%%%%%%%%%%%%%%%%%%%%%%%%%%%%%%%%%%%%%%%%%%%
From the above discussions, we can introduce an important lemmas for the $k$-slant helix in general form as follows:

\begin{lemma}\label{lm-13} If the Frenet frame of the spherical image of $\psi_{k}=$ indicatrix of the curve $\psi$ is $\{\t_k,\n_k,\b_k\}$, then we have Frenet formula:
\begin{equation}\label{u40}
 \left[
   \begin{array}{c}
     \t^{\,'}_k(s_k)\\
     \n^{\,'}_k(s_k)\\
     \b^{\,'}_k(s_k)\\
   \end{array}
 \right]=\left[
           \begin{array}{ccc}
             0 & \kappa_k & 0 \\
             -\kappa_k & 0 & \tau_k \\
             0 & -\tau_k & 0 \\
           \end{array}
         \right]\left[
   \begin{array}{c}
     \t_k(s_k)\\
     \n_k(s_k)\\
     \b_k(s_k)\\
   \end{array}
 \right],
 \end{equation}
 where
\begin{equation}\label{u41}
\t_k=\psi_{k+1},\,\,\,\,\,\n_k=\psi_{k+2},\,\,\,\,\,
\b_k=\frac{\psi_{k+1}\times\psi_{k+2}}{\|\psi_{k+1}\times\psi_{k+2}\|},
\end{equation}
and
\begin{equation}\label{u42}
\left\{
  \begin{array}{ll}
    s_k=\int\kappa(s)\sqrt{1+\sigma_0^2(s)}\sqrt{1+\sigma_1^2(s)}\,...\sqrt{1+\sigma_{k-1}^2(s)}\,ds,\\
    \kappa_k=\sqrt{1+\sigma_{k-1}^2},\\
    \tau_k=\sigma_k\sqrt{1+\sigma_{k-1}^2},
  \end{array}
\right.
\end{equation}
where
\begin{equation}\label{u43}
\sigma_k=\frac{\sigma'_{k-1}}{\kappa(s)\sqrt{1+\sigma_0^2(s)}\sqrt{1+\sigma_1^2(s)}\,...\,\Big(1+\sigma_{k-1}^2(s)\Big)^{3/2}},
\end{equation}
$s_k$ is the natural representation of the spherical image of $\psi_k$ indicatrix of the curve $\psi$ and $\kappa_k$ and $\tau_k$ are the curvature and torsion of this curve.
\end{lemma}

From the the above lemma we have $\frac{\tau_k}{\kappa_k}=\sigma_k$, which leads the following lemma:
\begin{lemma}\label{lm-131} Let $\psi:I\rightarrow\e^3$ be a $k$-slant helix. The spherical image of $\psi_{k}$ indicatrix of the curve $\psi$ is a spherical helix.
\end{lemma}

\begin{lemma}\label{lm-14} Let $\psi:I\rightarrow\e^3$ be a curve that is parameterized by arclength with intrinsic equations $\kappa(s)\neq0$ and $\tau(s)\neq0$. The curve $\psi$ is a $k$-slant helix (the vector $\psi_{k+1}$ makes a constant angle, $\phi$, with a fixed straight line in the space) if and only if the function \begin{equation}\label{u44}
\sigma_k=\cot[\phi].
\end{equation}
\end{lemma}

\begin{lemma}\label{lm-15} Let $\psi:I \rightarrow\e^3$ be a curve that is parameterized by arclength with intrinsic equations $\kappa(s)\neq 0$ and $\tau(s)\neq 0$. The curve $\psi$ is a $k$-slant helix if and only if the vector $\b_k=\frac{\psi_{k+1}\times\psi_{k+2}}{\|\psi_{k+1}\times\psi_{k+2}\|}$ makes a constant angle with fixed direction.
\end{lemma}

%%%%%%%%%%%%%%%%%%%%%%%%%%%%%%%%%%%%%%%


\begin{thebibliography}{99}
%%%%%%%%%%%%%%%%%%%%%%%%%%%%%%%%%%%%%%%

\bibitem{ali1} A.T. Ali: {\it Inclined curves in the Euclidean 5-space $\e^5$}, J. Adv. Res. in Pure Math. {\bf 1} (2009), 15--22.

\bibitem{ali2} A.T. Ali: {\it Determination of the position vector of general helices from intrinsic equations in $\e^3$},
Preprint 2009: arXiv:0904.0301v1 [math.DG].

\bibitem{ali3} A.T. Ali: {\it Position vectors of slant helices in Euclidean 3-space},
Preprint 2009: arXiv:0907.0750v1 [math.DG].

\bibitem{barros} M. Barros: {\it General helices and a theorem of Lancret}, Proc. Amer. Math. Soc. {\bf 125} (1997), 1503--1509.

\bibitem{bukcu} B. Bukcu and M.K. Karacan: {\it The slant helices according to Bishop frame}, Int. J. Comput. Math. Sci. {\bf 3} (2009), 67--70.

\bibitem{camci} C. Camci, K. Ilarslan, L. Kula and H.H. Hacisalihoglu: {\it Harmonic curvatures and generalized helices in $\e^n$}, Chaos, Solitons and Fractals {\bf 40} (2009), 2590--2596.

\bibitem{gluck} H. Gluck: {\it Higher curvatures of curves in Euclidean space},
Amer. Math. Monthly {\bf 73} (1996), 699--704.

\bibitem{hacis} H.H. Hacisalihoglu: Differential Geometry, Ankara University, Faculty of Science Press, 2000.

\bibitem{ilarslan} K. Ilarslan and O. Boyacioglu: {\it Position vectors of a spacelike W-cuerve in Minkowski space $\e_1^3$},
Bull. Korean Math. Soc. {\bf{44}} (2007), 429--438.

\bibitem{izumi} S. Izumiya and N. Takeuchi: {\it New special curves and developable surfaces},
Turk. J. Math. {\bf{28}} (2004), 531--537.

\bibitem{koend} J. Koenderink: Solid Shape, MIT Press, Cambridge, MA, 1990.

\bibitem{kuhn} W. Kuhnel: Differential Geometry: Curves - Surfaces - Manifolds. Wiesdaden: Braunchweig; 1999.

\bibitem{kula1} L. Kula and Y. Yayli: {\it On slant helix and its spherical indicatrix},
Appl. Math. Comp. {\bf{169}} (2005), 600--607.

\bibitem{kula2} L. Kula, N. Ekmekci, Y. Yayli and K. Ilarslan: {\it Characterizations of slant helices in Euclidean 3-space},
Tur. J. Math. {\bf{33}} (2009), 1--13.

\bibitem{mont2} J. Monterde: {\it Curves with constant curvature ratios}, Bulletin of Mexican Mathematic Society, 3a serie vol. {\bf 13} (20007), 177--186.

\bibitem{mont1} J. Monterde: {\it Salkowski curves revisted: A family of curves with constant curvature and non-constant torsion}, Comput. Aided Geomet. Design {\bf 26} (2009), 271--278.

\bibitem{salkow} E. Salkowski: {\it Zur transformation von raumkurven}, Mathematische Annalen {\bf 66} (1909), 517--557.

\bibitem{scofield} P.D. Scofield: {\it Curves of constant precession},
Amer. Math. Monthly {\bf 102} (1995), 531--537.

\bibitem{struik} D.J. Struik: Lectures in Classical Differential Geometry, Addison,-Wesley, Reading, MA, 1961.

\bibitem{turgut} M. Turgut and S. Yilmaz:
{\it Contributions to classical differential geometry of the curves in $\e^3$},
Scientia Magna {\bf 4} (2008),  5--9.

\end{thebibliography}
\end{document}